\documentclass[12pt]{article}
\usepackage{amsmath,amsthm,amssymb,amscd}
\usepackage{latexsym}
 \theoremstyle{definition}
 
 \theoremstyle{remark}

 \numberwithin{equation}{section}

\title
{Ricci flow on open  4-manifolds with positive isotropic curvature}

\author{ Hong Huang$$
  \\
\small School of Mathematical Sciences, Key Laboratory of Mathematics and Complex Systems,\\
\small Beijing Normal University \\
\small Beijing 100875,
P.R. China\\
\small {\em E-mail address:} {\bf hhuang@bnu.edu.cn}\\
}
\date{}
\begin{document}
\maketitle
\begin{abstract}
 In this note we prove the following result: Let $X$ be a
complete, connected 4-manifold with uniformly positive isotropic
curvature, with bounded geometry and with no essential incompressible space form. Then $X$ is diffeomorphic to
$\mathbb{S}^4$, or $\mathbb{RP}^4$, or $\mathbb{S}^3\times \mathbb{S}^1$, or $\mathbb{S}^3\widetilde{\times} \mathbb{S}^1$,
or a possibly infinite connected sum of them.
 This
extends   work of  Hamilton and Chen-Zhu to the noncompact case. The
proof uses Ricci flow with surgery on complete 4-manifolds, and is
inspired by recent work of Bessi$\grave{e}$res, Besson and
Maillot.


\end{abstract}
\maketitle

\section {Introduction}

In a recent paper [BBM] Bessi$\grave{e}$res, Besson and Maillot
classified complete 3-manifolds with uniformly positive scalar
curvature and with bounded geometry using a variant of
Hamilton-Perelman's Ricci flow with surgery. Inspired by
their work we try to classify complete 4-manifolds with uniformly
positive isotropic curvature, with bounded geometry and with no essential incompressible space form. More
precisely we will show

\hspace*{0.4cm}

{\bf Theorem }1.1. \ \  Let $X$ be a
complete, connected 4-manifold with uniformly positive isotropic
curvature, with bounded geometry and with no essential incompressible space form. Then $X$ is diffeomorphic to
$\mathbb{S}^4$, or $\mathbb{RP}^4$, or $\mathbb{S}^3\times \mathbb{S}^1$, or $\mathbb{S}^3\widetilde{\times} \mathbb{S}^1$, or a possibly infinite connected sum of them.

\hspace*{0.4cm}

 \noindent (Here, $\mathbb{S}^3\widetilde{\times} \mathbb{S}^1$ is the
only unorientable $\mathbb{S}^3$ bundle over $\mathbb{S}^1$. The
notion of  a (possibly infinite) connected sum will be given later
in this section; cf. [BBM].   By [MW] it is easy to see that the
converse is also true:
 Any 4-manifold as in the conclusion of the theorem has no essential incompressible space form, and admits a complete metric
  with uniformly positive isotropic curvature and with bounded
geometry.)

\noindent This extends  work of Hamilton [H5]  and Chen-Zhu [CZ2] to the noncompact
case.

Recall ([MM]) that a Riemannian manifold $M$ is said to have
positive isotropic curvature (PIC) if for all points $p \in M$ and
all orthonormal 4-frames $\{e_1,e_2,e_3,e_4\}\subset T_pM$ the
curvature tensor satisfies
\begin{equation*}R_{1313}+R_{1414}+R_{2323}+R_{2424}>2R_{1234}.
\end{equation*}

 Now  we consider in particular a 4-dimensional manifold $X$. If we decompose the bundle $\Lambda^2TX$
into the direct sum of its self-dual and anti-self-dual parts
\begin{equation*}
\Lambda^2TX=\Lambda_+^2TX \oplus \Lambda_-^2TX,
\end{equation*}
 then the curvature operator can be decomposed as
\begin{equation*}
\mathcal{R}=\left(
  \begin{array}{cc}
    A & B \\
    B^T & C \\
  \end{array}
\right),
\end{equation*}
 where $A=W_++\frac{R}{12}$, $C=W_-+\frac{R}{12}$, (here
$W_+ $ and $W_-$ are the self-dual part and the anti-self-dual part
of the Weyl curvature respectively,) and $B$ gives the trace free
part of the Ricci tensor.  Denote the eigenvalues of the matrices
$A, C$ and $\sqrt{BB^T}$ by $a_1\leq a_2 \leq a_3$, $c_1 \leq c_2
\leq c_3$ and $b_1 \leq b_2 \leq b_3$ respectively. It is easy to
see (cf. Hamilton [H5]) that  for a Riemannian 4-manifold the
condition of positive isotropic curvature is equivalent to the
condition $a_1+a_2>0$ and $ c_1+c_2>0$. A  Riemannian 4-manifold $X$
is said to have uniformly positive isotropic curvature if there is a
positive constant $c$ such that $a_1+a_2\geq c$ and $c_1+c_2\geq c$
everywhere.

As in [H5], an incompressible space form  in a 4-manifold $X$ is a 3-dimensional submanifold $Y$ diffeomorphic to $\mathbb{S}^3/\Gamma$
(where $\Gamma$ is a finite, fixed point free subgroup of isometries of $\mathbb{S}^3$) such that $\pi_1(Y)$ injects into $\pi_1(X)$.
The space form is called essential unless $\Gamma=1$, or $\Gamma=\mathbb{Z}_2$ and the normal bundle is non-orientable.
Also recall that a complete Riemannian manifold is said to have
bounded geometry if the sectional curvature is bounded (in both
sides) and the injectivity radius is bounded away from zero.

Now we explain the notion of (possibly infinite) connected sum,
following [BBM]. Let  $\mathcal{X}$ be a class of closed 4-manifolds. A
4-manifold $X$ is said to be a connected sum of members of
$\mathcal{X}$ if there exists a locally finite  graph $G$ and a map
$v\mapsto X_v$ which associates to each vertex of $G$ a copy of some
manifold in $\mathcal{X}$, such that by removing from each $X_v$ as
many open 4-balls as vertices incident to $v$ and gluing the thus
punctured $X_v$'s to each other along the edges of $G$ using
diffeomorphisms of the boundary 3-spheres, one obtains a 4-manifold
diffeomorphic to $X$.

Hamilton [H5] first used the Ricci flow with surgery to study
compact 4-manifolds with positive isotropic curvature and with no
essential incompressible space-form. (As Perelman [P2] pointed out,
[H5] contains some unjustified statements. See also [CZ2].) Later in
a breakthrough [P1], [P2] Perelman introduced some important new
ideas for the analysis of the Ricci flow, and devised a somewhat
different surgery procedure for it: one of the differences lies in
that Hamilton does surgery before curvature blows up, while Perelman
does surgery  exactly when curvature blows up. (For more details,
variants and/or alternatives of Perelman's arguments, see for
examples [BBB$^+$], [CaZ], [KL], [MT] and [Z].) Using Perelman's
ideas Chen-Zhu [CZ2] gave a complete proof of Hamilton's main
theorem in [H5]. Recently Chen-Tang-Zhu [CTZ] completely classified
all compact 4-manifolds (and 4-orbifolds with  isolated
singularities) with positive isotropic curvature using Ricci flow
with surgery on orbifolds. ( Note that in [CZ2] and [CTZ] the
surgeries are done exactly when curvature blows up as in [P2].)

Our proof of  Theorem 1.1 uses a 4-dimensional analogue of a version
of surgery constructed by Bessi$\grave{e}$res, Besson and Maillot
([BBM]) in 3-dimension;  see also [BBB$^+$].  Their  surgery
procedure is closer to that of Hamilton in the sense that they do
surgery before the curvature blows up; on the other hand, they also
use crucial ideas from Perelman [P1], [P2]. However, I adopt a
somewhat different approach from that in [BBM] to prove the
existence of $(r, \delta, \kappa)$-surgical solution with initial
data a complete 4-manifold with uniformly positive isotropic
curvature, with bounded geometry and with no essential
incompressible space form, see Theorem 3.4. Note that Perelman's
proof of [P2, Proposition 5.1] uses the openness (w.r.t. time)
property of canonical neighborhood assumption. In noncompact case it
is not clear whether it is still true. It turns out that a weak
openness (w.r.t. time) property of canonical neighborhood assumption
holds in our noncompact situation; see Claim 1 in the proof of
Proposition 3.6. We also need a slightly more general form of the
persistence of almost standard cap (in the phrase of [BBB$^+$]), see
Proposition 3.1, which corresponds to [P2, Lemma 4.5]. With these
tools in hand, we can adapt the original proof in [P2] and [CZ2] to
our noncompact case. Our approach can be adapted to treat more
general cases than that is considered in this note. Actually, I have
used the method in this note to deal with complete 4-orbifolds with
uniformly positive isotropic curvature, see [Hu1] and [Hu2]. (Those
two papers were written before this note was, and the main results
of this note are special case of those two papers, but I think maybe
it is worth to write down the details of this more simple case,
since in this case one needs not to worry about the additional
complexity in the orbifold case, and the main idea is clearer.) I
benefit much from [BBB$^+$], [BBM], [H5], [P1], [P2] and [CZ2]. In
particular, many definitions and proofs in this note are adapted
from [BBB$^+$], [BBM], [P2] and [CZ2].

In Section 2 we give some definitions and preliminary results, and
in Section 3,  we construct $(r, \delta,
\kappa)$-surgical solution with initial data a complete 4-manifold
with uniformly positive isotropic curvature, with bounded geometry
and with no essential incompressible space form, then Theorem 1.1
follows quickly. In Appendix A we collect some technical results on
gluing $\varepsilon$-necks, and finally in Appendix B we give a
version of bounded curvature at bounded distance for our surgical
solution, following [P1], [P2]. In most cases we will follow the
notations and conventions in [BBB$^+$] and [BBM].

\section{Surgical solutions on open 4-manifolds with uniformly PIC}

 Let $(X,g_0)$ be a complete 4-manifold with  $|Rm|\leq K$. Consider the Ricci
flow ([H1])
 \begin{equation}\frac{\partial g}{\partial t}=-2  Ric, \hspace{8mm}  g|_{t=0}=g_0.
\end{equation}
  \noindent By  Shi [S], (2.1) has a short time
solution with complete time slice and with bounded curvature. By Chen-Zhu ([CZ1]) this solution
is unique (in the category of complete solutions with bounded curvature).

 Now we assume  that the 4-manifold  $(X,g_0)$ has uniformly positive
isotropic curvature. Then we can easily generalize Hamilton's
pinching estimates in [H5] to our situation, which plays a similar role
in the category of 4-manifolds with uniformly positive isotropic
curvature as the Hamilton-Ivey pinching estimate does in the
category of 3-manifolds.

\hspace *{0.4cm}

  {\bf Lemma 2.1. }  (cf. Hamilton [H5])\ \ Let  $(X,g_0)$ be a complete 4-manifold with
  uniformly positive isotropic curvature ($a_1+a_2\geq c$, $c_1+c_2\geq c$) and with bounded curvature ($|Rm|\leq K$).  Then there exist positive
  constants $\varrho, \Psi, L, P, S<+\infty$ depending only on the initial
  metric  (through $c, K$), such that the complete solution to the Ricci flow (2.1) with bounded
  curvature satisfies
 \begin{equation}
 \begin{aligned}
& a_1+\varrho>0, \hspace{8mm} c_1+\varrho>0, \\
 & \max\{a_3,b_3,c_3\}\leq \Psi(a_1+\varrho), \hspace{8mm}  \max\{a_3,b_3,c_3\}\leq
  \Psi(c_1+\varrho), \\
& \frac{b_3}{\sqrt{(a_1+\varrho)(c_1+\varrho)}}\leq 1+\frac{L
e^{Pt}}{\max \{\ln \sqrt{(a_1+\varrho)(c_1+\varrho)},S\}}
\end{aligned}
\end{equation}

\noindent at all points and times.

\hspace *{0.4cm}

{\bf Proof}\ \  Note that Hamilton's maximum principle for Ricci
flow [H2] holds in the case of complete manifolds with bounded
curvature (see e.g. [CCG$^+$08, Chapter 12]). Then by inspecting
Hamilton's original proof in [H5,Section B] we see that the lemma is
true.

\hspace *{0.4cm}

Since the 4-manifolds we consider have uniformly positive isotropic
curvature, and in particular, have uniformly positive scalar
curvature, the Ricci flow (2.1) will blow up in finite time. Using
Lemma 2.1, we see that any blow-up limit (if it exists) coming from a solution as in Lemma 2.1 satisfies
the following restricted isotropic curvature pinching condition
\begin{equation} a_3\leq \Psi a_1, \hspace{8mm} c_3\leq \Psi c_1, \hspace{8mm} b_3^2\leq
a_1c_1,
\end{equation}
and in particular, has nonnegative curvature operator.

\hspace *{0.4cm}

Following [H5], [BBB$^+$] and [BBM], we will do surgery before the
curvature blows up.
 Roughly speaking, the surgery procedure is: start with (2.1),
at certain time before and near the first  time when the curvature will blow up, cutoff
necks in the manifold where the curvature is large,  glue back
 caps, and remove some components with known
topology to reduce the large curvature; continue the flow until one
comes near the next time when the curvature will  blow up, then do surgery as before, and
continue $\cdot\cdot\cdot$.

 Now we will adapt some definitions from
[BBM].

\hspace *{0.4cm}

{\bf Definition }([BBM]) \ \  Given an interval $I\subset
\mathbb{R}$, an evolving Riemannian manifold is a pair $(X(t),g(t))$
($t \in I$), where $X(t)$ is a (possibly empty or disconnected)
manifold  and $g(t)$ is a Riemannian metric on $X(t)$.  We say that
it is piecewise $C^1$-smooth if there exists a discrete subset $J$
of $ I$,
 such that the
following conditions are satisfied:

i. On each connected component of $I\setminus J$, $t \mapsto X(t)$
is constant (in topology), and $t \mapsto g(t)$ is $C^1$-smooth;

ii. For each $t_0\in J$, $X(t)=X(t_0)$ for any $t< t_0$ sufficiently
close to $t_0$, and $t\mapsto g(t)$ is left continuous at $t_0$;

iii. For each $t_0 \in J\setminus$  $\{$sup $I\}$, $t\mapsto
(X(t),g(t))$ has a right limit at $t_0$, denoted by
$(X_+(t_0),g_+(t_0))$.

\hspace *{0.4cm}

As in [BBM], a time $t\in I$ is regular if $t$ has a neighborhood in
$I$ where $X(\cdot)$ is constant and $g(\cdot)$ is $C^1$-smooth.
Otherwise it is singular.  We also denote by $f_{max}$ and $f_{min}$
the supremum and infimum of a function $f$, respectively, as in
[BBM].

\hspace *{0.4cm}

  {\bf Definition  }(Compare [BBM])\ \ A piecewise $C^1$-smooth
  evolving Riemannian 4-manifold $\{(X(t), g(t))\}_{t \in I }$  with uniformly positive isotropic curvature, with bounded curvature and with
  no essential incompressible space form  is called a
  surgical solution to the Ricci flow if it has the following
  properties.

  i. The equation $\frac{\partial g}{\partial t}=-2$ Ric is satisfied
  at all regular times;

  ii.  For each singular time $t$ one has $(a_1+a_2)_{min}(g_+(t))\geq
  (a_1+a_2)_{min}(g(t))$, $(c_1+c_2)_{min}(g_+(t))\geq
  (c_1+c_2)_{min}(g(t))$, and $R_{min}(g_+(t))\geq
  R_{min}(g(t))$;

  iii. For each singular time $t$ there is a locally finite collection
  $\mathcal{S}$ of disjoint, embedded $\mathbb{S}^3$'s in $X(t)$, and a manifold $X'$ such that

  (a) $X'$ is obtained from $X(t)\setminus \mathcal{S}$ by
  gluing back $B^4$'s  (closed 4-balls),

 (b) $X_+(t)$ is a union of some connected components of $X'$ and
 $g_+(t)=g(t)$ on $X_+(t)\cap X(t)$, and

(c) Each component of $X'\setminus X_+(t)$ is diffeomorphic to
$\mathbb{S}^4$, or $\mathbb{RP}^4$, or $\mathbb{RP}^4\sharp
\mathbb{RP}^4$, or $\mathbb{S}^3\times \mathbb{S}^1$, or
$\mathbb{S}^3\widetilde{\times} \mathbb{S}^1$, or $\mathbb{R}^4$, or
$\mathbb{RP}^4\setminus B^4$, or $\mathbb{S}^3\times \mathbb{R}$.

\hspace *{0.4cm}

{\bf Lemma 2.2.}  Any complete surgical solution with $a_1+a_2\geq
c$, $c_1+c_2\geq c$ and starting at $t=0$ must become extinct at
some time $T < \frac{1}{2c}$.

\hspace *{0.4cm}

{\bf Proof} From the evolution equation
\begin{equation}
\frac{\partial
R}{\partial t}=\triangle R+2|\text{Ric}|^2
\end{equation}
for the scalar curvature under Ricci flow, the maximum principle and
the definition above, any complete surgical solution with
$a_1+a_2\geq c$, $c_1+c_2\geq c$ must become extinct at some time
$T\leq \frac{2}{R_{min}(0)}< \frac{1}{2c}$.

\hspace *{0.4cm}

 Let $\{(X(t), g(t))\}_{t\in I}$ be a surgical solution
and $t_0\in I$. As in [BBM], if $t_0$ is singular, we set
$X_{reg}(t_0):= X(t_0)\cap X_+(t_0)$, and  $X_{sing}(t_0):=
X(t_0)\setminus X_{reg}(t_0)$. If $t_0$ is regular,
$X_{reg}(t_0)=X(t_0)$ and $X_{sing}(t_0)=\emptyset$. Let $t_0\in
[a,b]\subset I$ be a time, and $Y$ be a subset of $X(t_0)$ such that
for every $t \in [a,b)$, we have $Y\subset X_{reg}(t)$. Then as in
[BBM], we say the set $Y\times [a,b]$ is unscathed.

In [H5] Hamilton devised a quantitative metric surgery procedure;
later Perelman [P2] gave a somewhat different version, and in
particular, he had the crucial notion of ``canonical neighborhood''.
To describe it we need some more notions such as $\varepsilon$-neck,
$\varepsilon$-cap and strong $\varepsilon$-neck as given in [P2], [BBM],
[CZ2].

 Let $(X,g)$ be a Riemannian 4-manifold, and $x_0\in M$. An open neighborhood  $N\subset X$ of $x_0$  is an
$\varepsilon$-neck centered at $x_0$ if there is a diffeomorphism $\psi:
\mathbb{S}^3
  \times \mathbb{I} \rightarrow N$ such that the pulled back metric
  $\psi^*g$, scaling with some factor, is $\varepsilon$-close (in
  $C^{[\varepsilon^{-1}]}$ topology) to the standard metric $\mathbb{S}^3
  \times \mathbb{I}$ with scalar curvature 1 and
  $\mathbb{I}=(-\varepsilon^{-1},\varepsilon^{-1})$, and such that $x_0\in \psi(
\mathbb{S}^3 \times \{0\})$.

  An open subset
  $U$ is an $\varepsilon$-cap centered at $x_0$ if  $U$ is  the union of two sets $V$, $W$ such that $x_0\in$ Int $V$, $V$ is diffeomorphic to $B^4$  or
  $\mathbb{RP}^4\setminus $(Int $B^4$), $\overline{W}\cap V=\partial V$, and $W$ is an $\varepsilon$-neck.

    Let $(X(t), g(t))$ be an evolving Riemannian 4-manifold, and $(x_0, t_0)$ be a space-time point. An open subset $N\subset X(t_0)$ is a  strong
  $\varepsilon$-neck centered at  $(x_0,t_0)$ if there is a number $Q>0$ such that  the
  set  $\{(x,t)| x \in N, t\in [t_0-Q^{-1},t_0]\}$ is unscathed,
  and  there
  is a diffeomorphim $\psi: \mathbb{S}^3
  \times \mathbb{I} \rightarrow N$ such that, the pulled back
  solution $\psi^*g(\cdot,\cdot)$ scaling with the factor $Q$ and shifting
  the time $t_0$ to 0, is $\varepsilon$-close (in $C^{[\varepsilon^{-1}]}$
  topology) to the subset $(\mathbb{S}^3
  \times \mathbb{I})\times [-1,0]$ of the evolving round cylinder $\mathbb{S}^3
  \times \mathbb{R}$, with scalar curvature one and length 2$\varepsilon^{-1}$ to $\mathbb{I}$ at time zero, and $x_0\in \psi(\mathbb{S}^3 \times \{0\})$.

Motivated  by the structure theorems of 4-dimensional ancient
$\kappa$-solution with restricted isotropic curvature pinching
([CZ2, Theorem 3.8]) and the standard solution ([CZ2, Corollary
A.2]), following [P2], [BBM], [CZ2], we introduce the notion of
canonical neighborhood.

\hspace *{0.4cm}

{\bf Definition}  Let $\varepsilon$ and $C$ be positive constants. A
point $(x,t)$ in a surgical solution to the Ricci flow is said to
have an $(\varepsilon,C)$-canonical neighborhood if it has  an open
neighborhood $U$, $\overline{B_t(x,\sigma)} \subset U\subset
B_t(x,2\sigma)$ with $C^{-1}R(x,t)^{-\frac{1}{2}}<\sigma
<CR(x,t)^{-\frac{1}{2}}$, which falls into one of the following
three types:

(a) $U$ is a strong $\varepsilon$-neck with center $(x,t)$,

(b) $U$ is an $\varepsilon$-cap with center $x$ for $g(t)$,

(c) at time $t$, $U$ is a compact 4-manifold with positive curvature
operator,

\noindent and moreover, the scalar curvature in $U$ at time $t$ is
between $C^{-1}R(x,t)$ and $CR(x,t)$, and satisfies the derivative
estimates
\begin{equation*}
|\nabla R|< C R^{\frac{3}{2}} \hspace*{8mm} and \hspace*{8mm}
|\frac{\partial R}{\partial t}|< C R^2,
\end{equation*}
and the volume estimate
\begin{equation*}
(CR(x,t))^{-2}< vol_t(U).
\end{equation*}

\hspace *{0.4cm}

{\bf Remark }  Note that by [CZ2, Proposition 3.4 and Theorem 3.8]) and [CZ2, Corollary
A.2], for every $\varepsilon>0$, there exists a positive constant
$C(\varepsilon)$ such that each point in any ancient
$\kappa$-solution with restricted isotropic curvature
pinching  or in the standard
solution has an $(\varepsilon,C(\varepsilon))$-canonical
neighborhood, except that for the standard solution, an
$\varepsilon$-neck may not be strong.

\hspace *{0.4cm}

We choose $\varepsilon_0>0$ such that $\varepsilon_0<10^{-4}$ and such
that when $\varepsilon \leq 2\varepsilon_0$, Lemma A.1 in Appendix A and the results in the paragraph following its proof  hold true. Let $\beta:=\beta(\varepsilon_0)$ be the constant given by
Lemma A.2 in Appendix A.  Define $C_0:=\max\{100\varepsilon_0^{-1},
2C(\beta\varepsilon_0/2)\}$, where $C(\cdot)$ is given in the Remark
above. Fix $c_0>0$. Let $\varrho_0, \Psi_0, L_0, P_0, S_0$ be the
constants given in Lemma 2.1 by setting $c=c_0$ and $K=1$.

Now we consider some a priori assumptions, which consist of the
pinching assumption and the canonical neighborhood assumption.

\hspace *{0.4cm}

{\bf Pinching assumption}: Let
   $\varrho_0$, $\Psi_0$, $L_0$, $P_0$, $S_0$  be positive constants as given above. A  surgical solution to the Ricci flow
  satisfies the pinching assumption (with pinching constants $\varrho_0,\Psi_0,L_0, P_0, S_0$) if there hold
\begin{equation}\begin{aligned}
& a_1+\varrho_0>0, \hspace{8mm} c_1+\varrho_0>0, \\
&  \max\{a_3,b_3,c_3\}\leq \Psi_0(a_1+\varrho_0), \hspace{8mm}
\max\{a_3,b_3,c_3\}\leq \Psi_0(c_1+\varrho_0), \\
&  and  \\
& \frac{b_3}{\sqrt{(a_1+\varrho_0)(c_1+\varrho_0)}}\leq 1+\frac{L_0
e^{P_0t}}{\max\{\ln \sqrt{(a_1+\varrho_0)(c_1+\varrho_0)},S_0\}}
\end{aligned}\end{equation}
at all points and  times.

\hspace *{0.4cm}

{\bf Canonical neighborhood assumption}:   Let $\varepsilon_0$
and $C_0$ be as given above.
 Let $r: [0,+\infty)\rightarrow (0,+\infty)$ be a non-increasing function. An evolving Riemannian 4-manifold $\{(X(t), g(t))\}_{t \in I}$
 satisfies the canonical neighborhood assumption  $(CN)_r$ if  any  space-time point $(x,t)$ with  $R(x,t)\geq
r^{-2}(t)$ has an  $(\varepsilon_0,C_0)$-canonical
neighborhood.

\hspace *{0.4cm}

Let $\{(X(t), g(t))\}_{t \in I}$ be an evolving Riemannian
4-manifold. Recall [P1] that given $\kappa>0$, $r>0$, $g(\cdot)$ is
$\kappa$-noncollapsed  at $(x,t)$ (where $t\geq r^2$, and $P(x, t,
r, -r^2)$ is unscathed)  on the scale $r$ if
\begin{equation*}
|Rm|\leq r^{-2} \hspace{2mm} on \hspace{2mm} P(x, t, r, -r^2)
\hspace{2mm} implies \hspace{2mm} volB(x,t, r)\geq \kappa r^4,
\end{equation*}
where $P(x, t, r, -\Delta t):=\{(x', t')| x'\in B(x, t, r), t'\in
[t-\Delta t, t]\}$.

Let $\kappa: I \rightarrow (0, +\infty)$ be a function. We say $\{(X(t), g(t))\}_{t \in I}$ has property $(NC)_\kappa$ if it
is $\kappa (t)$-noncollapsed at any space-time point $(x,t)$ on all scales $\leq
1$.

 \hspace *{0.4cm}

The following proposition is analogous to [BBM, Theorem 6.5] and
[BBB$^+$, Theorem 6.2.1].

\hspace *{0.4cm}

{\bf Proposition 2.3.}\ \ Fix $c_0>0$.
 For any  $r$, $\delta>0$,
there exist    $h \in (0, \delta r)$ and $D> 10$, such that if
$(X(\cdot),g(\cdot))$ is a complete surgical solution  with
uniformly positive isotropic curvature ($a_1+a_2\geq c_0$, $c_1+c_2
\geq c_0$), with bounded curvature and with no essential
incompressible space form, defined on an time interval $[a,b]$
($0\leq a < b < \frac{1}{2c_0}$) and satisfying the pinching
assumption and the canonical neighborhood assumption $(CN)_r$, then
the following holds:

 \noindent Let $t \in [a,b]$ and  $x,y, z \in X(t)$ such that $R(x,t) \leq
2/r^2$,  $R(y,t)=h^{-2}$ and $R(z,t)\geq D/h^2$. Assume there is a
curve $\gamma$ in $X(t)$ connecting $x$ to $z$ via $y$, such that
each point of $\gamma$ with scalar curvature in $[2C_0r^{-2},
C_0^{-1}Dh^{-2}]$ is the center of an $\varepsilon_0$-neck.  Then
$(y,t)$ is the center of a strong $\delta$-neck.

\hspace *{0.4cm}

{\bf Proof}\ \ We  follow closely the proof of [BBM, Theorem 6.5]
and [BBB$^+$, Theorem 6.2.1]. (Compare [P2, Lemma 4.3], [CZ2, Lemma
5.2].) We argue by contradiction. Otherwise, there exist $r,
\delta>0$,  sequences $h_k \rightarrow 0$, $D_k\rightarrow +\infty$,
a sequence of complete surgical solutions $(X_k(\cdot),g_k(\cdot))$
 with uniformly positive isotropic
curvature ($a_1+a_2\geq c_0$, $a_1+a_2 \geq c_0$), with bounded
curvature and with no essential incompressible space form,
satisfying the  pinching assumption (with constants $\varrho_0,
\Psi_0, L_0, P_0, S_0$) and $(CN)_r$, and sequences $0 \leq t_k<
\frac{1}{2c_0}$, $x_k, y_k, z_k \in X_k(t_k)$ with $R(x_k,t_k)\leq
2r^{-2}$, $R(y_k, t_k)=h_k^{-2}$ and $R(z_k,t_k)\geq D_kh_k^{-2}$,
and finally a sequence of curves $\gamma_k$ in $X_k(t_k)$ connecting
$x_k$ to $z_k$ via $y_k$, whose points of scalar curvature in
$[2C_0r^{-2},C_0^{-1}D_kh_k^{-2}]$ are  centers of
$\varepsilon_0$-necks, but  $y_k$  is not the center of a strong
$\delta$-neck.

 Consider  the rescaled solution $(\overline{X}_k(\cdot), \bar{g}_k(\cdot))$, where $\bar{g}_k(\cdot)=h_k^{-2}g_k(t_k+h_k^2t)$.
  By Theorem  B.1 in Appendix B (and the Remark after Theorem B.1),
 for any $\rho>0$, there exists $\Lambda(\rho)>0$ and
 $k_0(\rho)>0$  such that the  ball $({B}(\bar{y}_k, 0, \rho)), \bar{g}_k(0))$ has scalar curvature bounded above by  $ \Lambda(\rho) $  for  $k> k_0(\rho)$.
 (Here and below, we adopt the convention in [BBM] to put a bar on the points when the relevant geometric quantities are computed w.r.t. the metric $\bar{g}_k$.)
 Combined with the canonical neighborhood assumption, it implies that
 the  parabolic neighborhoods  ${P}(\bar{y}_k, 0, \rho, -\frac{1}{2\Lambda(\rho)})$ are unscathed,
 with scalar curvature bounded above by $2\Lambda(\rho)$ for all
 $k\geq k_1(\rho)>k_0(\rho)$. By the pinching assumption, we get a uniform control of the curvature operator there. Using a local version of Hamilton's compactness theorem  (see
[BBB$^+$, Theorem C.3.3]), we see that (a subsequence of)
$(\overline{X}_k(0), \bar{g}_k(0), \bar{y}_k)$ converges to some
complete noncompact Riemannian 4-manifold $(\overline{X}_\infty,
\bar{g}_\infty, \bar{y}_\infty)$. Clearly $\overline{X}_\infty$ must
be diffeomorphic to $\mathbb{S}^3 \times \mathbb{R}$, and satisfy
(2.3) (with $\Psi=\Psi_0)$. By Toponogov's theorem it is the metric
product of some metric on $\mathbb{S}^3$ with $\mathbb{R}$; moreover
the spherical factor of this product must be $2\varepsilon_0$-close
to the round metric on $\mathbf{S}^3$ with scalar curvature 1. By
the closeness of the sequence $(\overline{X}_k(0), \bar{g}_k(0),
\bar{y}_k))$ to this limit  and properties of strong necks, for any
$\rho>0$, there exists $k_2(\rho)\geq k_1(\rho)$, such that for any
$k\geq k_2(\rho)$ the parabolic neighborhoods $P(\bar{y}_k, 0, \rho,
-\frac{1}{2})$ are unscathed, and have   scalar curvature satisfying
$\frac{1}{2}\leq R \leq 2$. By the local
 compactness theorem  ([BBB$^+$, Theorem C.3.3]) again, it follows that $(\overline{X}_k, \bar{g}_k(\cdot),
(\bar{y}_k,0))$ subconverges to some complete Ricci flow
$\bar{g}_\infty(\cdot)$ on $\overline{X}_\infty$. This  flow is
defined on $(-\frac{1}{2}, 0]$, has $R\leq 2$, and still satisfies
(2.3).

Now set
\begin{equation*}
\begin{aligned}
\tau_0:=\sup\{\tau>0| \forall \rho>0, \exists C(\rho, \tau)>0, \exists k(\rho), \forall k\geq k(\rho), P(\bar{y}_k, 0, \rho, -\tau) \\
 is \hspace{2mm} unscathed \hspace{2mm} and \hspace{2mm} C(\rho, \tau)^{-1} \leq R \leq C(\rho, \tau)  \hspace{2mm} there\}.
 \end{aligned}
 \end{equation*}
 We have shown $\tau_0\geq \frac{1}{2}$. It turns out  that, as in Step 2 of the proof of [BBB$^+$, Theorem 6.2.1], using the canonical neighborhood assumption one can show $\tau_0=+\infty$.
This way we get an ancient solution which satisfies (2.3) and splits at the final time slice. By
[CZ2, Lemma 3.2] it must be  the standard flow on the round cylinder. This implies the point $(y_k, t_k)$ is  the center of a strong $\delta$-neck when $k$ is sufficiently large --a contradiction.

\hspace *{0.4cm}

 Now we describe more precisely Hamilton's surgery procedure [H5].  We
will follow [CZ2] closely. First we describe the model surgery on
the standard cylinder, and define the standard solution. Consider
the semi-infinite cylinder $N_0=(\mathbb{S}^3 \times (-\infty,4)$
with the standard metric $\bar{g}_0$ of scalar curvature 1. Let $f$
be a smooth nondecreasing convex function on $(-\infty, 4)$ defined
by
\begin{equation*}
  \begin{cases}
    f(z)=0, & {z\leq 0;} \\
    f(z)=w_0e^{-\frac{W_0}{z}}, &{z\in (0,3];} \\
    f(z) \ \   \text {is strictly convex} , & {z\in [3,3.9];} \\
    f(z)=-\frac{1}{2}\text {ln} (16-z^2), & {z\in [3.9,4).}
  \end{cases}
  \end{equation*}
(where $w_0$ and $W_0$ are universal  positive
constants given in Lemma 2.4 below). Replace the standard
metric  $\bar{g}_0$ on the subspace $\mathbb{S}^3 \times [0,4)$ in
$N_0$ by $e^{-2f}\bar{g}_0$. The resulting metric will induce  a
complete, smooth metric (denoted by) $\hat{g}_0$ on $\mathbb{R}^4$. We call
the complete Ricci flow $(\mathbb{R}^4,\hat{g}(\cdot))$ with initial
data $(\mathbb{R}^4,\hat{g}_0)$ and with bounded curvature in any compact
subinterval of $[0, \frac{3}{2})$ the  standard solution, which
exists on the time interval $[0, \frac{3}{2})$. Denote by $p_0$ the
tip of the standard solution, which is the fixed point of the
$SO(4)$-action on the initial metric $(\mathbb{R}^4,\hat{g}_0)$.
 Note that by [CZ2,
Appendix], there exists a  constant $\kappa_{st}>0$ such that the
standard solution is $\kappa_{st}$-noncollapsed on scales $\leq 1$.
We refer the reader to [CZ2, Appendix] for other properties of
4-dimensional standard solution.

 Then we describe a similar surgery procedure for the general case.
 Suppose we have a $\delta$-neck centered at $x_0$ in a Riemannian 4-manifold $(X,g)$. Sometimes we will call $R^{-\frac{1}{2}}(x_0)$ the  radius of this neck.
 Let $\Phi: \mathbb{S}^3 \times [-l, l]\rightarrow
V\subset N$ be Hamilton's parametrization; see Appendix A. Assume the center $x_0$ of the
$\delta$-neck has $\mathbb{R}$ coordinate $z=0$. The surgery is to
cut off the $\delta$-neck along the middle 3-sphere and glue back
two balls (caps) separately.
 We construct a new smooth metric on the glued back cap  (say on the left hand side) as
follows.
\begin{equation*}
 \tilde{g}= \begin{cases}
    g(t_0), & {z=0;}  \\
    e^{-2f}g(t_0), & {z\in [0,2];} \\
    \varphi e^{-2f}g(t_0)+(1-\varphi)e^{-2f}h^2\bar{g}_0, & {z\in [2,3];} \\
    e^{-2f}h^2\bar{g}_0, & {z\in [3,4].}
  \end{cases}
\end{equation*}
where $\varphi$ is a smooth bump function with
$\varphi=1$ for $z\leq 2$, and $\varphi=0$ for $z\geq 3$, $h=R^{-\frac{1}{2}}(x_0)$, and
$\bar{g}_0$ is as above. We also perform the same surgery procedure
on the right hand side with parameter $\bar{z}\in [0,4]$
($\bar{z}=8-z$).

 The following  lemma of Hamilton justifies the pinching assumption
of surgical solution.

\hspace *{0.4cm}

 {\bf Lemma  2.4} \ \ (Hamilton [H5,Theorem D3.1]; compare  [CZ2, Lemma 5.3]) There
exist universal positive constants $\delta_0$, $w_0$ and $W_0$, and a constant $h_0$ depends only on $c_0$, such
that given any surgical solution with uniformly positive isotropic curvature ($a_1+a_2\geq c_0$, $c_1+c_2\geq c_0$),
satisfying the pinching assumption, defined on $[a, t_0]$ ($0\leq a< t_0< \frac{1}{2c_0}$), if we perform Hamilton's  surgery as described above at a
$\delta$-neck (if it exists) of radius $h$ at time $t_0$ with
$\delta <\delta_0$ and $h \leq h_0$, then after the surgery, the
pinching assumption still holds at
all points at time $t_0$. Moreover, after the surgery, any metric ball of radius $\delta^{-\frac{1}{2}}h$ with center near the tip
(i.e. the origin of the attached cap) is, after scaling with the factor $h^{-2}$, $\delta^{\frac{1}{2}}$-close to the corresponding ball of $(\mathbb{R}^4, \hat{g}_0)$.

\hspace *{0.4cm}

Usually we will be given two non-increasing step functions $r, \delta: [0,+\infty)\rightarrow (0, +\infty)$ as surgery parameters.  Let $h(r,\delta), D(r,\delta)$ be the
associated parameter  as determined in Proposition 2.3, ($h$ is also
called the surgery scale,) and let $ \Theta:=2Dh^{-2}$ be the
curvature threshold for the surgery process ( as in [BBM]), that is, we will do surgery only when $R_{max}$ reaches $\Theta$.

\hspace *{0.4cm}

Now we adapt two more definitions from [BBM].

\hspace *{0.4cm}

 {\bf Definition} (compare [BBM])\ \ Fix surgery parameter functions $r$, $\delta$  and let $h$, $D$,
$\Theta=2Dh^{-2}$ be the associated cutoff parameters. Let
$(X(t),g(t))$ ($ t \in I\subset
[0,\frac{1}{2c_0})$) be an evolving Riemannian 4-manifold with
uniformly positive isotropic curvature ($a_1+a_2\geq c_0$, $c_1+c_2
\geq c_0$), with bounded curvature and with no essential
incompressible space form. Let $t_0 \in I$ and $(X_+,g_+)$ be a
(possibly empty) Riemannian 4-manifold. We say that $(X_+,g_+)$ is
obtained from $(X(\cdot),g(\cdot))$ by $(r,\delta)$-surgery at time
$t_0$ if the following conditions are satisfied:

i. $R_{max}(g(t_0))=\Theta(t_0)$, and there is a locally finite
collection
  $\mathcal{S}$ of disjoint embedded $\mathbb{S}^3$'s in $X(t_0)$ which are in the middle  of strong $\delta$-necks with radius equal to the surgery scale $h(t_0)$, such that
  $X_+$ is obtained from $X(t_0)$ by doing
  Hamilton's surgery as described above on these necks, and removing the components that are
  diffeomorphic to  $\mathbb{S}^4$, or $\mathbb{RP}^4$, or $\mathbb{RP}^4\sharp
\mathbb{RP}^4$, or $\mathbb{S}^3\times \mathbb{S}^1$, or
$\mathbb{S}^3\widetilde{\times} \mathbb{S}^1$, or $\mathbb{R}^4$, or
$\mathbb{RP}^4\setminus B^4$, or $\mathbb{S}^3\times \mathbb{R}$.

ii. If $X_+\neq \emptyset$, then $R_{max}(g_+)\leq \Theta(t_0)/2$.

\hspace *{0.4cm}

{\bf Definition} (cf. [BBM])\ \ Fix surgery parameter functions $r$, $\delta$
and let $h$, $D$, $\Theta=2Dh^{-2}$ be the associated cutoff
parameters. A surgical solution $(X(\cdot),g(\cdot))$ with uniformly
positive isotropic curvature ($a_1+a_2\geq c_0$, $c_1+c_2 \geq
c_0$), with bounded curvature and with no essential incompressible
space form, defined on some time interval $I\subset
[0,\frac{1}{2c_0})$ is an $(r,\delta)$-surgical solution  if it has
the following properties:

i. It satisfies the pinching assumption, and $R(x,t) \leq \Theta(t)$
for all $(x,t)$;

ii. At each singular time $t_0\in I$, $(X_+(t_0),g_+(t_0))$ is
obtained from $(X(\cdot),g(\cdot))$ by $(r,\delta)$-surgery at time
$t_0$; and

iii. Condition $(CN)_r$ holds.

\noindent Let $\kappa$ be a positive function (here, usually a nonincreasing step function). An $(r,\delta)$-surgical solution which also satisfies
Condition $(NC)_\kappa$  is called an
$(r,\delta,\kappa)$-surgical solution.

\hspace *{0.4cm}

The following lemma is analogous to [BBM, Lemma 5.9].

\hspace *{0.4cm}

{\bf Lemma 2.5}\ \ Suppose we have fixed two constants  $r, \delta>0$ as surgery parameters on an interval $[a,b)$. Let $(X(t),g(t))$ be an
$(r, \delta)$-surgical solution on $[a,b]$. Let $a \leq t_1 <
t_2 < b$ be two singular times (if they exist). Then $t_2-t_1$ is
bounded from below by a positive number depending only on
$r,\delta$.

\hspace *{0.4cm}

{\bf Proof}\ \ We may assume that there are no other singular times
between $t_1$ and $t_2$. Since $R_{max}(g_+(t_1))\leq \Theta/2$,
$R_{max}(g(t_2))=\Theta$, and $\Theta$ depends only on $r, \delta$,
the result follows by integrating  the curvature derivative estimate $|\frac{\partial R}{\partial t}|< C_0 R^2$ in the canonical neighborhood assumption (see [BBM, Lemma 5.9]).

\hspace *{0.4cm}

The following proposition is similar to [BBM, Theorem 7.4], and it
extends a result in [CZ2] to the noncompact case.

\hspace *{0.4cm}

{\bf Proposition 2.6}\ \ Let $\varepsilon\in (0, 2\varepsilon_0]$.
Let $(X,g)$ be a complete, connected 4-manifold. If each point of
$X$ is the center of an $\varepsilon$-neck or an $\varepsilon$-cap,
then $X$ is diffeomorphic to $\mathbb{S}^4$, or $\mathbb{RP}^4$, or
$\mathbb{RP}^4\sharp \mathbb{RP}^4$, or $\mathbb{S}^3\times
\mathbb{S}^1$, or $\mathbb{S}^3\widetilde{\times} \mathbb{S}^1$, or
$\mathbb{R}^4$, or $\mathbb{RP}^4\setminus B^4$, or
$\mathbb{S}^3\times \mathbb{R}$.

\hspace *{0.4cm}

{\bf Proof}. The result in the compact case has been shown in [CZ2].
So below we will assume that $X$ is not compact.

\hspace *{0.4cm}

{\bf Claim }  Let $\varepsilon\in (0, 2\varepsilon_0]$. Let $(X,g)$
be a complete, noncompact, connected 4-manifold. If each point of
$X$ is the center of an $\varepsilon$-neck, then $X$ is
diffeomorphic to
 $\mathbb{S}^3 \times \mathbb{R}$.

\hspace *{0.4cm}

Proof of Claim. Let $x_1$ be a point of $X$, and let $N_1$ be a
$\varepsilon$-neck centered at $x_1$,  given by some diffeomorphism
$\psi_1:\mathbb{S}^3 \times (-\varepsilon^{-1},
\varepsilon^{-1})\rightarrow N_1$. Consider Hamilton's canonical
parametrization $\Phi_1: \mathbb{S}^3 \times [-l_1,l_1]\rightarrow
V_1\subset N_1$ such that $V_1$ contains the portion
$\psi_1(\mathbb{S}^3 \times (-0.98\varepsilon^{-1},
0.98\varepsilon^{-1}))$ in $N_1$. (See Appendix A.)  Now choose a
point $x_2$ in $\Phi_1( \mathbb{S}^3 \times \{0.9l_1\})$, and let
$N_2$ be a $\varepsilon$-neck centered at $x_2$, given by some
diffeomorphism $\psi_2:\mathbb{S}^3 \times (-\varepsilon^{-1},
\varepsilon^{-1})\rightarrow N_2$. Again consider Hamilton's
canonical parametrization $\Phi_2: \mathbb{S}^3 \times
[-l_2,l_2]\rightarrow V_2\subset N_2$ such that $V_2$ contains the
portion $\psi_2(\mathbb{S}^3 \times (-0.98\varepsilon^{-1},
0.98\varepsilon^{-1}))$ in $N_2$.   Then  by [H5, Theorem C2.4]  we
have Hamilton's canonical parametrization $\Phi: \mathbb{S}^3 \times
[-l,l]\rightarrow  V_1\cup V_2$, and for all $\alpha \in [-l_1,
l_1]$ and all $\beta \in [-l_2, l_2]$, $\Phi_1(\mathbb{S}^3 \times
\{\alpha\})$  is isotopic to $\Phi_2(\mathbb{S}^3 \times
\{\beta\})$. (See also
 Appendix A.) Then we go on, choose $x_3$, $N_3$, $\Phi_3$,
$\cdot\cdot\cdot.$ This way the desired result follows.

\hspace *{0.4cm}

Now consider the case  that $X$ contains at least one
$\varepsilon$-cap. In this case, since we are assuming $X$ is
noncompact, $X$ contains only one cap. Then arguing as above, one
see that $X$ is diffeomorphic to a cap. So in this case $X$ is diffeomorphic to $\mathbb{R}^4$  or $\mathbb{RP}^4\setminus B^4$.

\hspace *{0.4cm}

The following proposition is  analogous to [BBM, Proposition A].

\hspace *{0.4cm}

{\bf Proposition 2.7}\ \ Fix $c_0>0$. There exists a positive
constant $\tilde{\delta}$ (depending only on $c_0>0$) with the
following property: Let $r, \delta$ be surgery parameters, let
$\{(X(t), g(t))\}_{t\in (a,b]}$ ( $0<a<b<\frac{1}{2c_0}$) be an $(r,
\delta)$-surgical solution  with uniformly positive isotropic
curvature ($a_1+a_2\geq c_0$, $c_1+c_2 \geq c_0$), with bounded
curvature, and with no essential incompressible space form. Suppose
that $\delta\leq \tilde{\delta}$, and $R_{max}(b)=\Theta(b)$.
 Then there exists a Riemannian manifold
$(X_+,g_+)$ which is obtained from $(X(\cdot),g(\cdot))$ by
$(r,\delta)$-surgery at time $b$, such that

i.  $g_+$ satisfies the pinching assumption at time $b$;

ii. $(a_1+a_2)_{min}(g_+(b))\geq
  (a_1+a_2)_{min}(g(b))$, $(c_1+c_2)_{min}(g_+(b))\geq
  (c_1+c_2)_{min}(g(b))$, and $R_{min}(g_+(b))\geq R_{min}(g(b))$;

iii. $X_+$ has no essential incompressible space form.

 \hspace *{0.4cm}

{\bf Proof}\ \  Let  $\delta_0$ and
$h_0$ be as  given in Lemma 2.4. Set
$\tilde{\delta}=\frac{1}{2} \min \{c_0^{\frac{1}{2}}h_0,
\delta_0\}$.

For the proof of i. and ii. we will follow that of [BBM, Proposition
A]. Let $\mathcal{G}$ (resp. $\mathcal{O}$, resp. $\mathcal{R}$) be
the set of points of $X(b)$ of scalar curvature less than $2r^{-2}$
(resp. $\in [2r^{-2}, \Theta(b)/2)$, resp. $\geq \Theta(b)/2$). The idea
is to consider a maximal collection $\{N_i\}$ of pairwise disjoint
cutoff necks in $X(b)$, whose existence is guaranteed by Zorn's
Lemma . (Here, following [BBM], a cutoff neck is a strong
$\delta$-neck centered at some point $(x,b)$ with
$R(x,b)=h^{-2}$.) It is easy to see that such a collection is locally
finite by a volume argument.

\hspace *{0.4cm}

{\bf Claim 1} Any connected component of $X(b)\setminus \cup_iN_i$
is contained either in $\mathcal{G}\cup \mathcal{O}$ or in
$\mathcal{R} \cup \mathcal{O}$.

\hspace *{0.4cm}

Proof of Claim 1. We argue by contradiction. Otherwise there is some
component $W$ of $X(b)\setminus \cup_i N_i$ containing at least one
point $x \in \mathcal{G}$ and one point $z \in \mathcal{R}$. Choose
a minimizing geodesic path $\gamma$ in $W$ connecting $x$ with $z$.
In the following Claim 2, we will show each point of $\gamma$ with
scalar curvature in $[2C_0r^{-2}, C_0^{-1}Dh^{-2}]$ is the center of
an $\varepsilon_0$-neck.  Then we can apply Proposition 2.3 to
conclude that there exists some point $y\in \gamma$ with
$R(y,b)=h^{-2}$ which is the center of a strong $\delta$-neck. This
will contradict the maximality of $\{N_i\}$.

\hspace *{0.4cm}

{\bf Claim 2}  Each point of such $\gamma$ with scalar curvature in
$[2C_0r^{-2}, C_0^{-1}Dh^{-2}]$ is the center of an
$\varepsilon_0$-neck.

\hspace *{0.4cm}

Proof of Claim 2. The proof is a minor modification of that of the
second claim in Lemma 7.7 of [BBM]. Let $y\in \gamma$ be such a
point. Then $y$ is the center of an $(\varepsilon_0,C_0)$-canonical
neighborhood  $U$. Clearly $U$ cannot be a closed manifold by the
curvature assumptions. We will show $U$ cannot be an
$(\varepsilon_0,C_0)$-cap either. Otherwise $U=N\cup C$, where $N$
is an $\varepsilon_0$-neck, $N\cap C=\emptyset$, $\overline{N}\cap
C=\partial C$ and $y\in$ Int $C$. Let $\psi: \mathbb{S}^3 \times
(-\varepsilon_0^{-1}, \varepsilon_0^{-1}) \rightarrow N$ be the
diffeomorphism which defines the neck $N$. We use Hamilton's method
to give a canonical parametrization $\Phi: \mathbb{S}^3 \times
[-l,l]\rightarrow V\subset N$ such that $V$ contains the portion
$\psi(\mathbb{S}^3 \times (-0.98\varepsilon_0^{-1},
0.98\varepsilon_0^{-1}))$  (cf. Lemma A.1 in  Appendix A). Let
$S=\Phi( \mathbb{S}^3 \times \{0\})$. We rescale the metric such
that the scalar curvature of $N$ is close to 1. Clearly $\gamma$ is
not minimizing in $U$, since if $x'$ (resp. $z'$) is an intersection
of $\gamma$ with $S$ between $x$ and $y$ (resp. $y$ and $z$), then
$d(x',z')\ll d(x',y)+d(y,z')$. The geodesic segment (in $U$)
$[x'z']$ is not contained in $W$ by the minimality of $\gamma$ in
$W$. So $ [x'z']\cap
\partial W \neq \emptyset$. By definition of $W$, the corresponding component of
$\partial W$ is a boundary component, denoted by $S_i^+$, of some
cutoff neck $N_i$. Then $d(S_i^+,S)< $diam$(S)$ since $ [x'z']\cap
S_i^+ \neq \emptyset$. We use Hamilton's method to give a canonical
parametrization $\Phi': \mathbb{S}^3 \times [-l', l']\rightarrow
V'\subset N_i$ such that one of the ends of $V'$, denoted by $\partial_+V'$,  is at the rescaled
distance  $< 0.03\varepsilon_0^{-1}$ from  the end $S_i^+$ of $N_i$.
Pick a point $p'$ in $V'$ which is at rescaled distance $0.2\varepsilon_0^{-1}$ from  $\partial_+V'$.
 Then $d(p', S)\leq d(p',\partial_+V')+ d(\partial_+V',S_i^+)+d(S_i^+,S)<0.03\varepsilon_0^{-1}+0.2\varepsilon_0^{-1}+$diam$(S)<0.3\varepsilon_0^{-1}$.
 Then it follows from the discussion after Lemma A.1  that  the embedded $\mathbb{S}^3$ in the neck structure of $V'$ which  contains $p'$ is isotopic to $S$ in $N$.
 It follows that $\gamma \cap N_i \neq \emptyset$, which is impossible by the definition of $W$.

\vspace *{0.4cm}

 Then we do Hamilton's surgery along these $N_i$'s, and obtain an
manifold $(X',g_+)$. The components of $X'$ consist of two types:
Either they have curvature $\leq \Theta(b)/2$, or they are covered by
canonical neighborhoods, whose diffeomorphism types are identified
with the help of Proposition 2.6, and will be thrown away.  We
denote the resulting manifold by  $(X_+,g_+)$.  By Lemma 2.4 and our
choice of $\tilde{\delta}$ it satisfies the pinching assumption.
Clearly ii) is also satisfied.

Now we  show it satisfies iii. also.   We will adapt an argument in
[CZ2] to the noncompact case. We argue by contradiction. Suppose
$X_+$ has  an essential incompressible space form $Y\approx
\mathbb{S}^3/\Gamma$, where $\Gamma$ is a finite, fixed point free
subgroup of isometries of $\mathbb{S}^3$. After an isotopy, we may
assume the intersection of $Y$ with the union of all surgery caps is
empty. Then $Y$ may be seen as a submanifold in $X(b)$ also. Below
we will show $Y$ is also an essential incompressible space form in
$X(b)$,  which contradicts to our assumption on $X(\cdot)$ and
completes the proof.

\vspace *{0.4cm}

{\bf Claim 3}  \hspace*{4mm}    $Y$ is also an essential incompressible space form in
$X(b)$.

\vspace *{0.4cm}

Proof of Claim 3. We argue by contradiction. Suppose  $Y$ is not an
essential incompressible space form in $X(b)$.

 Case 1. $Y$ is compressible in $X(b)$. Then we can pick a loop $\gamma
\subset Y$ representing a nontrivial element in the kernel of
$i_*:\pi_1(Y)\rightarrow \pi_1(X(b))$, where $i$ is the inclusion
map. So there is a map $f: D^2\rightarrow X(b)$ with $f(\partial
D^2)=\gamma$.  Since $f(D^2)$ is compact and the collection of our
cutoff necks is locally finite, $f(D^2)$ will intersects only a
finite number  of 3-spheres which lie in the middle of cutoff necks.
Denote these 3-spheres by $S_1$, $S_2$, $\cdot\cdot\cdot, S_m$. We
perturb they slightly so that they meet $f(D^2)$ transversely in a
finite number of simple closed curves. By using an innermost circle
argument we may assume (after modifying $f$ suitably) that the
enclosed disks in $D^2$ of all the circles in the preimage (of these
intersection curves) are disjoint; denote these circles by $C_1,
C_2, \cdot\cdot\cdot, C_l$, and the enclosed 2-disks by $D_1, D_2,
\cdot\cdot\cdot, D_l$. Each  $f(C_j)$ bounds a homotopical 2-disk in
$S_1\cup S_2\cup \cdot\cdot\cdot \cup S_m$, since each $S_k$ is a
topological 3-sphere. So after a further modification of $f$ we may
assume that  $f(D_1 \cup D_2 \cup \cdot\cdot\cdot \cup D_l)$ is
contained in $S_1\cup S_2\cup \cdot\cdot\cdot \cup S_m$. On the
other hand, since $D^2 \setminus ( D_1 \cup D_2 \cup \cdot\cdot\cdot
\cup D_l)$ is connected, $f(\partial D^2)=\gamma \subset Y$, we see
that $f(D^2 \setminus ( D_1 \cup D_2 \cup \cdot\cdot\cdot \cup
D_l))\subset X_+$.  So $\gamma$  bounds a homotopical disk in $X_+$.
This contradicts to the choice of $Y$.

Case 2.  $Y$ is incompressible in $X(b)$, but not essential. If
$\Gamma=\{1\}$ then $Y$ cannot be essential in $X_+$. So we may
assume  $\Gamma=\mathbb{Z}_2$  and the normal bundle of $Y$ in
$X(b)$ is non-orientable. But the normal bundle of $Y$ in $X(b)$ is
the same as in $X_+$. So $Y$ again cannot be essential in $X_+$. A
contradiction.

\section{Existence of $(r, \delta, \kappa)$-surgical solutions}

As in [BBM], if $(X(\cdot), g(\cdot))$ is a piecewise $C^1$ evolving
manifold defined on some interval $I\subset \mathbf{R}$ and
$[a,b]\subset I$, the restriction of $g$ to $[a,b]$, still denoted
by $g(\cdot)$, is the evolving manifold
\begin{equation*}
 t\mapsto \begin{cases}
     (X_+(a), g_+(a)), &       t=a, \\
    (X(t), g(t)), &          t\in (a,b]. \\
 \end{cases}
 \end{equation*}

The following proposition is analogous to [P2, Lemma 4.5] and [BBM, Theorem 8.1], which is
one of the key technical results in the process of constructing  $(r,
\delta, \kappa)$-surgical solutions; compare [BBB$^+$, Theorem
8.1.2], [CaZ, Lemma 7.3.6], [KL, Lemma 74.1],
[MT, Proposition 16.5] and [Z, Lemma 9.1.1], see also the
formulation in the proof of [CZ2, Lemma 5.5]. We state it in a slightly more general form, which is applicable to our situation.

\hspace *{0.4cm}

{\bf Proposition 3.1}\ \ Fix $c_0>0$. For all $A>0, \theta \in (0,\frac{3}{2})$
and $\hat{r}>0$, there exists
$\hat{\delta}=\hat{\delta}(A,\theta,\hat{r})>0$ with the following
property. Let $r(\cdot)\geq \hat{r}$, $\delta(\cdot)
\leq \hat{\delta}$ be two positive step functions on $[a,b)$ ($0\leq a< b<\frac{1}{2c_0}$), and let $(X(\cdot),g(\cdot))$ be a surgical solution  with
uniformly positive isotropic curvature ($a_1+a_2\geq c_0$,
$c_1+c_2 \geq c_0$), with bounded curvature and
with no essential incompressible space form, defined on $[a,b]$, such that it satisfies the pinching assumption on $[a,b]$,
that $R(x,t)\leq \Theta(r(t), \delta(t))$ for all space-time points with $t \in [a,b)$, that at any singular time $t_0\in [a,b)$,
$(X_+(t_0), g_+(t_0))$ is obtained from $(X(\cdot), g(\cdot))$ by  $(r, \delta)$-surgery, and that any point
$(x,t)$ ($t\in [a,b)$) with $R(x,t) \geq (\frac{r(t)}{2})^{-2}$ has a $(2\varepsilon_0, 2C_0)$-canonical neighborhood, Let $t_0\in [a,b)$ be a singular time. Consider
the restriction of $(X(\cdot),g(\cdot))$ to $[t_0,b]$. Let $p\in
X_+(t_0)$ be the tip of some surgery cap of scale $h(t_0)$, and let $t_1\leq$ min
$\{b,t_0+\theta h^2(t_0)\}$ be maximal (subject to this inequality) such that
$P(p,t_0,Ah(t_0),t_1-t_0)$ is unscathed. Then the following holds:

i. The parabolic neighborhood $P(p,t_0,Ah(t_0),t_1-t_0)$ is, after
scaling with factor $h^{-2}(t_0)$ and shifting time $t_0$ to zero,
$A^{-1}$-close to $P(p_0,0,A,(t_1-t_0)h^{-2}(t_0))$ (where
$p_0$ is the tip of the  standard solution);

ii. If $t_1<$ min $\{b,t_0+\theta h^2(t_0)\}$, then
$B(p,t_0,Ah(t_0))\subset X_{sing}(t_1)$ disappears at time $t_1$.

\hspace *{0.4cm}

We will follow the proof of [BBB$^+$, Theorem 8.1.2] and [BBM,
Theorem 8.1].

 Let $\mathcal{M}_0=(\mathbb{R}^4,\hat{g}(\cdot))$ be
the standard solution, and  $0<T_0<\frac{3}{2}$.

\hspace *{0.4cm}

The following result is from [BBB$^+$], where the proof uses Chen-Zhu's uniqueness theorem ([CZ1]).

 \hspace *{0.4cm}

{\bf Lemma 3.2} ( [BBB$^+$, Theorem 8.1.3])  For all $A, \Lambda>0$,
there exists $\rho=\rho(\mathcal{M}_0,A,\Lambda)>A$ with the
following property.  Let $U$ be an open subset of $\mathbb{R}^4$ and
$T\in (0,T_0]$. Let $g(\cdot)$ be a Ricci flow defined on $U\times
[0,T]$, such that the ball $B(p_0,0,\rho)\subset U$ is relatively
compact. Assume that

i. $||Rm(g(\cdot))||_{0,U\times [0,T],g(\cdot)}\leq \Lambda$,

ii. $g(0)$ is $\rho^{-1}$-close to $\hat{g}(0)$ on $B(p_0,0,\rho)$.

Then $g(\cdot)$ is $A^{-1}$-close to $\hat{g}(\cdot)$ on
$B(p_0,0,A)\times [0,T]$.

\hspace *{0.4cm}

Here, $||Rm(g(\cdot))||_{k,U\times [0,T],g(\cdot)}:= \sup_{U \times
[0,T]} \{|\nabla^iRm_{g(t)}|_{g(t)} |0\leq i \leq k\}$.

\hspace *{0.4cm}

{\bf Corollary 3.3} (Compare [BBM, Corollary 8.3])   Let $A>0$.
There exists $\rho=\rho(\mathcal{M}_0,A)
>A$ with the following property. Let $\{(X(t),g(t))\}_{t\in
[0,T]}$ ($T \leq T_0$) be a surgical solution
 with
uniformly positive isotropic curvature, with bounded curvature and
with no essential incompressible space form. Assume that

i. $(X(\cdot),g(\cdot))$ is a parabolic rescaling of some surgical solution which satisfies the pinching assumption,

ii.  $|\frac{\partial R}{\partial t}|\leq 2C_0 R^2$ at any space-time
point $(x,t)$ with $R(x,t)\geq 1$.

\noindent Let $p\in X(0)$ and $t\in (0,T]$ be such that

iii. $B(p,0,\rho)$ is $\rho^{-1}$-close to $B(p_0,0,\rho)$,

iv. $P(p,0,\rho,t)$ is unscathed.

\noindent Then $P(p,0,A,t)$ is $A^{-1}$-close to $P(p_0,0,A,t)$.

\hspace *{0.4cm}

{\bf Proof} The proof is similar to that of Corollaries 8.2.2 and
8.2.4 in [BBB$^+$].

\hspace *{0.4cm}

Using  Corollary 3.3, one can easily adapt  the arguments in the proof of [BBB$^+$, Theorem
8.1.2] and [BBM, Theorem 8.1] to prove Proposition 3.1.

\hspace *{0.4cm}

The following theorem is analogous to [P2, Proposition 5.1] and [BBM, Theorems 5.5 and 5.6]. We state it in a form similar to [MT, Theorem 15.9].

\hspace *{0.4cm}

{\bf Theorem 3.4}  \ \ Given  $c_0$,
$v_0>0$, there are surgery parameter sequences
\begin{equation*}
\mathbf{K}=\{\kappa_i\}_{i=1}^\infty, \hspace{2mm}
\Delta=\{\delta_i\}_{i=1}^\infty, \hspace{2mm}
\mathbf{r}=\{r_i\}_{i=1}^\infty
\end{equation*}
such that the following holds. Let $r(t)=r_i$ and
$\bar{\delta}(t)=\delta_i$  on $[(i-1)2^{-5}, i\cdot2^{-5})$,
$i=1, 2, \cdot\cdot\cdot$. Suppose that $\delta:
[0,\infty)\rightarrow (0,\infty)$ is a non-increasing step function with
$\delta(t) \leq \bar{\delta}(t)$. Then the following holds: Suppose
that  we have a surgical solution $(X(\cdot), g(\cdot))$ with
uniformly positive isotropic curvature, with bounded curvature and with no essential incompressible space form, defined on $[0,T]$ (for some $T< \infty$), which
satisfies the following conditions:

(1) the initial data $(X(0),g(0))$ is  a complete 4-manifold  with
uniformly positive isotropic curvature ($a_1+a_2\geq c_0$, $c_1+c_2\geq c_0$), with $|Rm|\leq 1$, with no essential incompressible space form, and with vol $B(x,1)\geq v_0$ at any point $x$,

(2) the solution satisfies the pinching assumption, and $R(x,t)\leq \Theta(r(t),\delta(t))$ for all space-time points,

(3) it has only a finite  number of singular times  such that at
each singular time $t_0\in (0,T)$, $(X_+(t), g_+(t))$ is obtained
from $(X(\cdot),g(\cdot))$ by $(r,\delta)$-surgery at time
$t_0$, and

(4) on each time interval $[(i-1)2^{-5}, i\cdot2^{-5}]\cap
[0,T]$ the solution satisfies $(CN)_{r_i}$ and $(NC)_{\kappa_i}$.

\noindent Then there is an extension of $(X(\cdot), g(\cdot))$ to a
surgical solution  defined
for  $0\leq t\leq T'$ (where $T'<\frac{1}{2c_0}$ is the extinction
time) and satisfying the above four conditions with $T$ replaced by
$T'$.

\hspace *{0.4cm}

To prove the theorem above, I will adapt the arguments in Perelman [P2] and Chen-Zhu [CZ2] to the  noncompact case; compare [BBM].

\hspace *{0.4cm}

The following lemma guarantees the non-collapsing under a weak form of the canonical
neighborhood assumption,  and is analogous to [P2, Lemma 5.2], [BBM, Proposition C], and [CZ2, Lemma 5.5]. We state it in a
form close to [ KL, Lemma 79.12].

 \hspace *{0.4cm}

{\bf Lemma 3.5}\ \  Fix $c_0>0$. Suppose $0< r_- \leq \varepsilon_0 $, $\kappa_->0$,
and $0< E_-<
 E< \frac{1}{2c_0}$. Then there exists
 $\kappa_+=\kappa_+(r_-,\kappa_-,E_-,E)>0$, such that for any  $r_+$, $0< r_+ \leq r_-$, one can find
  $\delta'=\delta'( r_-,r_+,\kappa_-,E_-,E)>0 $,  with the following
 property.

 \noindent Suppose that $0\leq a<b<d < \frac{1}{2c}$, $b-a\geq E_-$, $d-a\leq E$. Let $r$ and $\delta
$ be two positive step functions on $[a,d)$ with $\varepsilon_0\geq
r\geq r_-$ on $[a,b)$,
 $\varepsilon_0\geq r\geq r_+$ on $[b,d)$  and $\delta \leq \delta'$ on $[a,d)$.
Let  $(X(\cdot),g(\cdot))$ be a surgical solution
 with
uniformly positive isotropic curvature ($a_1+a_2\geq c_0$, $c_1+c_2
\geq c_0$), with bounded curvature and with no essential
incompressible space form,  defined on the time interval $[a, d]$,
such that it satisfies the pinching assumption on $[a, d]$, that
$R(x,t)\leq \Theta(r(t),\delta(t))$ for all space-time points with
$t\in [a,d)$, that at any singular time $t_0\in [a,d)$, $(X_+(t_0),
g_+(t_0))$ is obtained from $(X(\cdot), g(\cdot))$ by $(r,
\delta)$-surgery, that the conditions $(CN)_{r}$   and
$(NC)_{\kappa_-}$ hold on $[a,b)$, and that any point $(x,t)$ ($t\in
[b, d)$) with $R(x,t) \geq (\frac{r(t)}{2})^{-2}$ has a
$(2\varepsilon_0, 2C_0)$-canonical neighborhood. Then
$(X(\cdot),g(\cdot))$ satisfies $(NC)_{\kappa_+}$ on $[b,d]$.

\hspace *{0.4cm}

\hspace *{0.4cm}

{\bf Proof} \ \ Using Proposition 3.1  and Perelman's reduced
volume, the proof of [CZ2,Lemma 5.2] can be adapted to our case
without essential changes.

\hspace *{0.4cm}

The following proposition justifies the canonical neighborhood
assumption needed.  We state it in a form similar to [MT, Proposition 17.1]. Compare  [P2, Section 5], [BBM, Proposition B] and  [CZ2, Proposition 5.4].

\hspace *{0.4cm}

{\bf Proposition 3.6} Given $c_0>0$.
 Suppose that for some $i\geq 1$ we have surgery parameter sequences
 $\tilde{\delta}\geq \delta_1\geq \delta_2 \geq \cdot\cdot\cdot\geq \delta_i >0$, $\varepsilon_0 \geq r_1 \geq\cdot\cdot\cdot \geq r_i>0$ and
 $\kappa_1 \geq \kappa_2 \geq\cdot\cdot\cdot \geq \kappa_i>0$, where $\tilde{\delta}$ is the constant given in Proposition 2.7.
 Then there are positive constants $r_{i+1}\leq r_i$ and $\delta_{i+1}\leq \min \{\delta_i, \delta'\}$,
 where $\delta'=\delta'(r_i, r_{i+1}, \kappa_i)$ is the constant given in Lemma 3.5 by setting $r_-=r_i$, $\kappa_-=\kappa_i$, $r_+=r_{i+1}$, $E_-=2^{-5}$ and $E= 2^{-4}$,
 such that the following holds.  Let $r(t)=r_j$ and
$\bar{\delta}(t)=\delta_j$  on $[(j-1)2^{-5}, j\cdot2^{-5})$,
$j=1, 2, \cdot\cdot\cdot, i+1$. Suppose that $\delta:
[0, (i+1)2^{-5})\rightarrow (0,\infty)$ is a non-increasing step function with
$\delta(t) \leq \bar{\delta}(t)$.  Let
$(X(\cdot),g(\cdot))$ be any surgical solution to Ricci flow with uniformly positive isotropic curvature
($a_1+a_2\geq c_0$, $c_1+c_2\geq c_0$), with bounded curvature and with no essential incompressible space form, defined on $[0, T]$ for some $T\in
(i\cdot2^{-5}, (i+1)2^{-5}]$, such that $R(x,t)\leq \Theta(r(t),\delta(t))$ for all space-time points with $t\in [0,T)$, that there are only a
finite number of singular times, and at each singular time $t_0\in
(0,T)$, $(X_+(t_0),g_+(t_0))$ is obtained from $(X(\cdot),g(\cdot))$
by $(r,\delta)$-surgery at time $t_0$. Suppose that the
restriction of the surgical solution to $[0, i\cdot2^{-5}]$
satisfies the four conditions given in Theorem 3.4. Suppose also that
$\delta(t) \leq \delta_{i+1}$ for all $t \in
[(i-1)2^{-5}, T)$. Then $(X(\cdot),g(\cdot))$ satisfies the
condition $(CN)_{r_{i+1}}$ on $[i\cdot2^{-5}, T]$.

\hspace *{0.4cm}

{\bf Proof} \ \  We argue by contradiction. Otherwise there exist
$r^\alpha \rightarrow 0$ as $\alpha \rightarrow \infty$, and for
each $\alpha$ a sequence  ${\delta}^{\alpha\beta}\rightarrow 0$ as
$\beta \rightarrow \infty$, such that the following holds. For each
$\alpha, \beta$ there is a surgical solution $(X^{\alpha\beta
}(\cdot), g^{\alpha\beta}(\cdot))$ to the Ricci flow defined for $0\leq t \leq
T_{\alpha\beta}$ with $i\cdot 2^{-5} < T_{\alpha\beta} \leq
(i+1)2^{-5}$, such that it satisfies the conditions of the
proposition w.r.t. these constants but not the conclusion.

\noindent

 By our assumption, Lemma 2.1 and Lemma 2.4, $(X^{\alpha\beta}(\cdot)
,g^{\alpha\beta}(\cdot))$ satisfies the pinching assumption on $[0,
T_{\alpha\beta}]$.  Note that by our assumption the scalar
curvature of $(X^{\alpha\beta }(t),g^{\alpha\beta}(t))$ on $[0,
i\cdot 2^{-5}]$ are uniformly bounded above by a constant
independent of $\alpha, \beta$ (but depending on $i$). So by
choosing $r^\alpha $ sufficiently small we may assume that the
condition $(CN)_{r^\alpha}$ holds on $[i\cdot 2^{-5}, i\cdot
2^{-5}+\theta]$ for some $\theta>0$. Also note that if for
$(X^{\alpha\beta}(\cdot),g^{\alpha\beta}(\cdot))$ the condition
$(CN)_{r^\alpha}$  holds in $[i\cdot 2^{-5}, \lambda)$ for some
$\lambda> i\cdot 2^{-5}$, then arguing as in the proof of [MT, Lemma
11.23] with some minor modifications, (note that the curvature of
$(X^{\alpha\beta}(t), g^{\alpha\beta}(t))$ are bounded on $[i\cdot
2^{-5}, \lambda)$,)  we see that any point $(x,t)$ ($t\in  [i\cdot
2^{-5}, \lambda]$) with $R(x,t) \geq (\frac{r^\alpha}{2})^{-2}$ has
a $(2\varepsilon_0, 2C_0)$-canonical neighborhood; cf. also
[BBB$^+$, Chapter 9]. This means that  the canonical neighborhood
condition has some sort of weak closeness (w.r.t. the time).

\noindent  The following Claim 1 may be seen as some sort of weak
openness (w.r.t. the time) property of the canonical neighborhood
condition.

 \hspace *{0.4cm}

{\bf Claim 1} \ \ Suppose  for   $(X^{\alpha\beta}(t),g^{\alpha\beta}(t))$ the condition  $(CN)_{r^\alpha}$
holds on $[i\cdot 2^{-5}, t_0]$ for some $i\cdot 2^{-5}< t_0
<T_{\alpha\beta}$. Then there exists $\tau>0$  (depending on
$\alpha, \beta$) such that  any point $(x,t)$ ($t\in  [i\cdot 2^{-5}, t_0+\tau]$) with $R(x,t) \geq (\frac{r^\alpha}{2})^{-2}$ has a $(2\varepsilon_0, 2C_0)$-canonical neighborhood.

 \hspace *{0.4cm}

{\bf Proof of Claim 1}\ \ We consider the following two cases.

Case i: $R_{max}(g^{\alpha\beta}(t_0))=\Theta=\Theta(r^\alpha,
\delta^{\alpha\beta})$ (the curvature threshold for the surgery
process). So  a surgery occurs  at $t_0$. Then
$R_{max}(g_+^{\alpha\beta}(t_0))\leq \Theta/2$, and the condition $(CN)_{r^\alpha}$  still holds in $(X_+^{\alpha\beta}(t_0),
g_+^{\alpha\beta}(t_0))$  (cf. for example the proof of [KL, Lemma
73.7]). Also note that the curvature derivatives of
$(X_+^{\alpha\beta}(t_0), g_+^{\alpha\beta}(t_0))$ are bounded. The
reason is as follows.
 If a point $x\in X_+^\alpha(t_0)$ lies
within a distance of $10\varepsilon_0^{-1}h$ from the added part
$X_+^{\alpha\beta}(t_0)\setminus X^{\alpha\beta}(t_0)$, (where
$h=h(r^\alpha, \delta^{\alpha\beta})$ is the surgery scale,) then it
lies in an $\varepsilon_0$-cap, and the curvature derivatives at
$(x,t_0)$ are bounded by the surgery construction. (Compare [MT,
Claim 16.6].) If a point $x\in X_+^{\alpha\beta}(t_0)$ lies at
distance greater than or roughly equal to $10\varepsilon_0^{-1}h$ from
$X_+^{\alpha\beta}(t_0)\setminus X^{\alpha\beta}(t_0)$, then it has
existed for a previous time interval since the set of singular times is a discrete subset of $\mathbb{R}$, and by Shi's local estimates [S] the
curvature derivatives at $(x,t_0)$ are bounded also. Then by Shi's
theorem with initial curvature derivative bounds ([LT,Theorem 11], see also [MT, Theorem
3.29]), the  curvature derivatives of
 $(X^{\alpha\beta}(t), g^{\alpha\beta}(t))$ are bounded when restricted to $[t_0,
 t'']$ for some $t''> t_0$ with $t''-t_0$ sufficiently small.

Case ii: $R_{max}(g^{\alpha\beta}(t_0))<\Theta=\Theta(r^{\alpha},
\delta^{\alpha\beta})$.
 By the smoothing property of Ricci flow there exists
$t''>t_0$ such that we have  $R_{max}(g^{\alpha\beta}(t))< \Theta$
on $[t_0, t'']$. If there exist singular times before $t_0$, let $t'$ be the last one. Otherwise let $t'=(i-1)2^{-5}+2^{-6}$. Then there are no surgeries in the time interval $(t', t'']$. By Shi's
theorem with  or without initial curvature derivative bounds  we see that the curvature derivatives of  $(X^{\alpha\beta}(t), g^{\alpha\beta}(t))$ are bounded on $[t_0,
 t'']$.

Arguing this way, we see the following holds:  Fix $0< \sigma< t_0$
and $t''> t_0$ with $t''-t_0$ sufficiently small as above. For any
$k\in \mathbb{N}$, the $k$-th covariant  derivative of the curvature
of
 $(X^{\alpha\beta}(t), g^{\alpha\beta}(t))$ are bounded when restricted to a subinterval $[t_1, t_2]$ of $[\sigma,
 t'']$ with $t_2-t_1$ sufficiently small. The bound depends  on $k$, $\alpha, \beta$, $\sigma$ and $t''$, but is independent of $t_1$ and $t_2$.  Then arguing as
in the proof of [H4, Lemma 2.4],   for any $k\in \mathbb{N}$, we
have sup$_{x\in
X^{\alpha\beta}(t_2)}\sum_{j=0}^k|\nabla_{g^{\alpha\beta}(t_1)}^j(g^{\alpha\beta}(t_2)-g^{\alpha\beta}(t_1))(x)|_{g^{\alpha\beta}(t_1)}\leq
c_k\cdot(t_2-t_1)\rightarrow 0$ as  $t_2\rightarrow t_1$, where
$c_k$ is a constant (depending also on $\alpha, \beta$, $\sigma$ and
$t''$, but not on $t_1$ or $t_2$). Moreover, from the curvature
derivative estimates  and the evolution equation (2.4), we know that
if $t-t_0>0$ is sufficiently small and $R(x,t)\geq
(r^{\alpha}/2)^{-2}$ for some $x\in X^{\alpha\beta}(t)$, then
$R(x,t_0)\geq (r^{\alpha})^{-2}$, and $(x, t_0)$ has an
$(\varepsilon_0, C_0)$-canonical neighborhood. Then Claim 1 follows
by combining the above estimates on metrics with the definition of
canonical neighborhood.

\hspace *{0.4cm}

From Claim 1, the weak closeness (w.r.t. the time) property of the
canonical neighborhood condition mentioned above and the choice of
our $(X^{\alpha\beta }(t),g^{\alpha\beta}(t))$ we have the following

\hspace *{0.4cm}

{\bf Claim 2} \ \ For each $(X^{\alpha\beta
}(t),g^{\alpha\beta}(t))$, there exists $t^{\alpha\beta}\in
(i\cdot 2^{-5}, T_{\alpha\beta}]$ such that
  the condition $(CN)_{r^\alpha}$  is violated at some
space-time point $(x^{\alpha\beta},t^{\alpha\beta})$, but any point $(x,t)$ ($t\in [i\cdot 2^{-5}, t^{\alpha\beta}]$) with
$R(x,t) \geq (\frac{r^\alpha}{2})^{-2}$ has a $(2\varepsilon_0, 2C_0)$-canonical neighborhood.

\hspace *{0.4cm}

 \noindent Note that we are not claiming that
$t^{\alpha\beta}$ is the first time that  the condition $(CN)_{r^\alpha}$ is violated.

Given $\alpha$, we may assume that for all
$\beta$,  $\delta^{\alpha\beta} \leq \delta'(r_i,r^\alpha/2,\kappa_i)$, the constant given  in Lemma 3.5 by setting  $r_-=r_i$, $r_+=r^{\alpha}/2$,
$\kappa_-=\kappa_i$, $E_-=2^{-5}$ and $E= 2^{-4}$. Then Claim 2 allows one to apply Lemma 3.5  to get uniform
$\kappa$-noncollapsing on all scales $\leq$ $1$ in $[0,
t^{\alpha\beta}]$ with  $\kappa=\kappa_+(r_i,\kappa_i)
>0$ independent of $\alpha, \beta$. (Note
that the $\kappa$-noncollapsed condition is  closed w.r.t. the time
(cf. [BBB$^+$, Lemma 4.1.4]).)

Now similarly as in [P2], [CZ2],  the  idea is roughly as follows: Let
$\bar{g}^{\alpha\beta}$ be the solutions obtained by rescaling
$g^{\alpha\beta}$ with factors
$R(x^{\alpha\beta},t^{\alpha\beta})$ and shifting the
times $t^{\alpha\beta}$ to 0. If for any $A>0$, $b>0$, the sets
$B(\bar{x}^{\alpha\beta}, 0, A) \times [-b,0]$ are
unscathed when $\alpha, \beta$ are sufficiently large,
 then with the help of  the  uniform $\kappa$-non-collapsing just obtained, the compactness theorem for Ricci flow, and Hamilton's Harnack estimate [H3], we can get an ancient
$\kappa$- solution with restricted isotropic curvature pinching as a
limit,
  which will lead to a
contradiction by the remark  (after the definition of canonical
neighborhood) in Section 2.  If there exist $A>0$, $b>0$ such that
there are arbitrarily large $\alpha, \beta$ with
$B(\bar{x}^{\alpha\beta}, 0, A) \times [-b,0]$ scathed, then using
Proposition 3.1 one can show that the solution will be close to the
standard solution (after suitable rescaling and time-shifting),
which again leads to a contradiction. (Compare also, for example,
the descriptions in [CaZ], [KL] and [MT].)

\noindent Actually one can argue similarly as in the proof of [CZ2,
Proposition 5.4] with some minor modifications. (One can also check
that the presentation in, for example,  [Z, Theorem 9.2.1] can be
adapted to our situation.)  I only indicate some of the  modifications (to arguments of [CZ2,
Proposition 5.4]) needed in our situation.

1. One can use Theorem B.1 in Appendix B to simplify or replace some
arguments in [CZ2, Proposition 5.4] (for example, Step 3 of the
proof there).

2. One may use the neck-strengthening lemma A.2 in  Appendix A to
replace the remark after Lemma 5.2 in [CZ2] used in the second
paragraph on p. 245 there.

\hspace *{0.4cm}

{\bf Proof of Theorem 3.4} (Compare Chapter 17, Section 2 in [MT].)
We will construct the desired  surgery parameter sequences inductively. By
Shi's work and the doubling time estimate (cf. [CLN, Lemma
6.1]), the initial condition that $(X, g_0)$ is complete
and has $|Rm|\leq 1$ guarantees a (smooth) complete  solution exists
on the time interval $[0, 2^{-5}]$, such that $|Rm|\leq 2$ on
this time interval. By Lemma 2.1, the pinching assumption is  satisfied on this time interval. Since we have the initial time curvature bound $|Rm|\leq 1$ and  the volume lower bound on the
initial time unit balls, by Perelman's no-local-collapsing theorem (cf. [KL, Theorem 26.2]), the solution satisfies $(NC)_{\kappa_1}$  for some
$\kappa_1>0$ on this time interval. By choosing $r_1\leq \varepsilon_0$, it satisfies
the condition $(CN)_{r_1}$
vacuously on this time interval. Pick any positive constant $\delta_1 \leq \tilde{\delta}$.
Now suppose we have constructed the surgery parameter sequences
\begin{equation*}
\mathbf{r}_i=\{r_1,\cdot\cdot\cdot, r_i\}, \hspace{2mm} \Delta_i=\{\delta_1,
\cdot\cdot\cdot, \delta_i\}, \hspace{2mm}  \mathbf{K}_i=\{\kappa_1,
\cdot\cdot\cdot, \kappa_i\},
\end{equation*}
with the desired property. We let $r_{i+1}$ and
$\delta_{i+1}$ be as in Proposition 3.6.  Then let
$\kappa_{i+1}=\kappa_+(r_i, \kappa_i)$ be as given in Lemma 3.5. Set
\begin{equation*}
\mathbf{r}_{i+1}=\{\mathbf{r}_i, r_{i+1}\},\hspace{2mm}
\Delta_{i+1}=\{\delta_1, \cdot\cdot\cdot, \delta_{i-1},
\delta_{i+1}, \delta_{i+1}\},\hspace{2mm}\mathbf{K}_{i+1}=\{\mathbf{K}_i,
\kappa_{i+1}\}.
\end{equation*}

\noindent Let $\delta:[0,(i+1)2^{-5}] \rightarrow (0,\infty)$ be
any non-increasing positive step function with $\delta \leq
\mathbf{\Delta}_{i+1}$.  Let $(X(\cdot),g(\cdot))$ be any surgical
solution to Ricci flow
 defined on
$[0, T]$ with $T\in [i\cdot 2^{-5}, (i+1)2^{-5})$ satisfying
the four conditions w.r.t.  $\delta$ and these sequences. We
want to extend this surgical solution to one defined on $[0,(i+1)2^{-5}]$ such that it still satisfies the
four conditions w.r.t. $\delta$,  $\mathbf{r}_{i+1}$,
$\Delta_{i+1}$, and $\mathbf{K}_{i+1}$.

If  $R_{max}(g(T))<\Theta=\Theta(r_{i+1}, \delta_{i+1})$, then we
can run the Ricci flow with initial data $(X(T), g(T))$ for a while
$[T, T+\theta]$ ($T+\theta \leq (i+1)2^{-5}$)  until
$R_{max}(T+\theta)$ reaches the threshold $\Theta$. By Lemma 2.1 the pinching assumption is satisfied.  Note that by
Proposition 3.6 the extended surgical flow still satisfies the
condition $(CN)_{r_{i+1}}$. (There are no
surgeries occurring in $[T, T+\theta)$, and the conditions of
Proposition 3.6 hold.)  Then by Lemma 3.5 it is
$\kappa_{i+1}$-noncollapsed.    So the extended flow satisfies
all the four conditions w.r.t. these constants.

So we may assume that $R_{max}(g(T))=\Theta=\Theta(r_{i+1},
\delta_{i+1})$.
  Using Proposition 2.7 we do $(r_{i+1}, \delta_{i+1})$-surgery and get
$(X_+(T), g_+(T))$. If $X_+(T)=\emptyset$, then we are done.  Assume $X_+(T)\neq \emptyset$, then we have
$R_{max}(g_+(T))\leq \frac{1}{2}\Theta$, $g_+(T)$ satisfies the pinching assumption at time $T$, and $X_+(T)$ has  no essential incompressible space form.
Then we run the Ricci flow with initial data $(X_+(T), g_+(T))$ for
a little while $[T, \lambda]$ for some  $\lambda> T$ keeping
$R_{max}(g(t))<\Theta$ on $(T,\lambda]$. By Lemma 2.1 the pinching assumption is satisfied. By Proposition 3.6
 the condition $(CN)_{r_{i+1}}$ still holds for the
extended flow. Then by Lemma 3.5 it is
$\kappa_{i+1}$-noncollapsed. So the extended flow satisfies
all the four conditions w.r.t. these constants.

In this way we can continue  the surgical solution to the time
interval $[0,(i+1)\cdot 2^{-5}]$ such that it always satisfies the four
conditions, since by Lemma 2.5  the
singular times cannot accumulate.

Finally note that by Lemma 2.2 the surgical solutions we have constructed must be extinct before the time $t=\frac{1}{2c_0}$.

\hspace *{0.4cm}

The following lemma  is essentially due to [BBM].

\hspace *{0.4cm}

{\bf Lemma 3.7} (see [BBM, Proposition 2.6])   Let $\mathcal{X}$ be a
class of closed 4-manifolds, and  $X$ be a 4-manifold. Suppose there
exists a finite sequence of 4-manifolds $X_0$, $X_1$,
$\cdot\cdot\cdot$, $X_p$ such that $X_0=X$, $X_p=\emptyset$, and for
each $i$ ($1\leq i \leq p$), $X_i$ is obtained from $X_{i-1}$ by
cutting off along a locally finite collection of pairwise disjoint,
embedded 3-spheres, gluing back $B^4$'s, and removing some
components that are connected sums of members of $\mathcal{X}$. Then
each component of $X$ is an connected sum of members of
$\mathcal{X}$.

\hspace *{0.4cm}

{\bf Proof} The proof  of [BBM, Proposition 2.6] applies to
4-dimensional case.

\hspace *{0.4cm}

 Note that each manifold  appeared in the list of our Proposition 2.6 is a  (possibly infinite)  connected sum of members of
 $\mathcal{X}=\{\mathbb{S}^4, \mathbb{RP}^4, \mathbb{S}^3\times \mathbb{S}^1, \mathbb{S}^3\widetilde{\times} \mathbb{S}^1\}$.  Then Theorem 1.1 follows from Theorem 3.4 and Lemma 3.7.

\section* {Appendix A}

The following lemma is essentially due to Hamilton [H5].

\hspace *{0.4cm}

{\bf Lemma A.1} \ \  Let  $\varepsilon
>0$ be sufficiently small. Suppose that $N$ is an  $\varepsilon$-neck centered
 at $x$, with a diffeomorphism  $\psi: \mathbb{S}^3 \times (-\varepsilon^{-1}, \varepsilon^{-1}) \rightarrow
N$,
 in a Riemannian
4-manifold  $(X,g)$. Then we have Hamilton's canonical
parametrization $\Phi: \mathbb{S}^3 \times [-l, l] \rightarrow
V\subset N$, such that  $V$ contains the portion $\psi(\mathbb{S}^3
\times (-0.98\varepsilon^{-1}, 0.98\varepsilon^{-1}))$ in $N$.

\hspace *{0.4cm}

{\bf  Proof}\ \  If $\varepsilon >0$ is sufficiently small, by the
inverse function theorem
 every point ( in the $\varepsilon$-neck $N$) at  distance (w.r.t the rescaled metric $R(x)g$) at least 0.01$\varepsilon^{-1}$ from
 the ends lies on a unique constant mean curvature hypersurface. Then we can choose harmonic parametrizations of the spheres,
 choose the height function, and straighten out the parametrization by rotations of the horizontal spheres to obtain  Hamilton's canonical parametrization
$\Phi: \mathbb{S}^3 \times [-l,l]\rightarrow  V\subset N$ such that
$V$ contains the portion $\psi(\mathbb{S}^3 \times
(-0.98\varepsilon^{-1}, 0.98\varepsilon^{-1}))$ in $N$. (Cf. the
proof of [H5, Theorem C2.2].)

\hspace *{0.4cm}

Now  suppose that $N_i$ is an  $\varepsilon$-neck centered
 at $x_i$, with a diffeomorphism  $\psi_i: \mathbb{S}^3 \times (-\varepsilon^{-1}, \varepsilon^{-1}) \rightarrow
N_i$, $i=1,2,$
 in a Riemannian
4-manifold $(X,g)$.
 Let $\pi_i: N_i \rightarrow (-\varepsilon^{-1}, \varepsilon^{-1})$
be the composition of $\psi_i^{-1}$ with the projection of
$\mathbb{S}^3 \times (-\varepsilon^{-1}, \varepsilon^{-1})$ onto its
second factor. Assume that $N_1 \cap N_2$ contains a point $y$ with
$-0.9 \varepsilon^{-1} \leq \pi_i(y) \leq 0.9\varepsilon^{-1}$
($i=1,2$). Then by the above lemma we have Hamilton's canonical
parametrization $\Phi_i: \mathbb{S}^3 \times [-l_i, l_i] \rightarrow
V_i\subset N_i$, such that  $V_i$ contains the portion
$\psi_i(\mathbb{S}^3 \times (-0.98\varepsilon^{-1},
0.98\varepsilon^{-1}))$ in $N_i$. If $\varepsilon$ is sufficiently
small, we can use [H5, Theorem C2.4] to get Hamilton's canonical
parametrization $\Phi: \mathbb{S}^3 \times [-l,l]\rightarrow V_1\cup
V_2$ and diffeomorphisms $F_1$ and $F_2$ of the cylinders, such that
$\Phi_1=\Phi\circ F_1$ and $\Phi_2=\Phi\circ F_2$. $F_1$ and $F_2$
are in fact isometries in the standard metrics on the cylinders by
[H5, Lemma C2.1]. Moreover we know that for all $\alpha \in [-l_1,
l_1]$ and all $\beta \in [-l_2, l_2]$, $\Phi_1(\mathbb{S}^3 \times
\{\alpha\})$  is isotopic to $\Phi_2(\mathbb{S}^3 \times
\{\beta\})$.

\hspace *{0.4cm}

The following lemma is essentially due to [BBB$^+$] and [BBM].

Let $K_{st}$ be the superemum of the sectional curvatures of the
(4-dimensional)  standard solution on $[0,4/3]$.

\hspace *{0.4cm}

 {\bf Lemma A.2}
 (see [BBB$^+$,Lemma 4.3.5] and [BBM,Lemma 4.11]) For any $\varepsilon \in (0, 10^{-4})$
 there exists $\beta=\beta (\varepsilon)\in (0,1)$ with the
 following property.

  Let $a, b$ be real numbers satisfying $a< b <0$ and $|b| \leq
  \frac{3}{4}$, let $(X(\cdot), g(\cdot))$ be a surgical solution (with no essential incompressible space form)
  defined on $(a,0]$, and $x \in X(b)$ be a point such that:

  i. $R(x,b)=1$;

  ii. $(x,b)$ is the center of a strong $\beta \varepsilon$-neck;

  iii. $P(x, b, (\beta \varepsilon)^{-1}, |b|)$ is unscathed and
  satisfies $|Rm|\leq 2K_{st}$.

 \noindent  Then $(x,0)$ is center of a strong $\varepsilon$-neck.

\hspace *{0.4cm}

{\bf Proof} The proof is almost identical to that of [BBB$^+$,Lemma 4.3.5].

\section* {Appendix B}

Bounded curvature at bounded distance is one of the key ideas in
Perelman [P1], [P2]; a 4-dimensional version appeared in [CZ2]. The
following 4-dimensional version is very close to [MT, Theorem 10.2],
[BBB$^+$, Theorem 6.1.1] and [BBM, Theorem 6.4].

\hspace *{0.4cm}

{\bf Theorem B.1} For each $c, \varrho, \Psi, L, P, S, A, C >0$ and
each $\varepsilon \in (0,2\varepsilon_0]$, there exists $Q=Q(c,
\varrho, \Psi,L, P, S,A, \varepsilon, C)>0$ and $\Lambda=\Lambda(c,
\varrho, \Psi, L, P,S, A, \varepsilon, C)>0$ with the following
property. Let $I=[a, b]$ ($0 \leq a < b < \frac{1}{2c})$ and
$\{(X(t), g(t))\}_{t\in I}$ be a surgical solution with uniformly
positive isotropic curvature ($a_1+a_2\geq c$, $c_1+c_2\geq c$),
with bounded curvature, with no essential incompressible space form
and satisfying the pinching condition  (2.2) (with constants
$\varrho, \Psi, L, P, S$). Let $(x_0, t_0)$ be a space-time point
such that:

1. $R(x_0, t_0)\geq Q$;

2. For each point $y\in B(x_0, t_0, AR(x_0, t_0)^{-1/2})$, if $R(y,
t_0)\geq 4R(x_0, t_0)$, then $(y, t_0)$ has an $(\varepsilon,
C)$-canonical neighborhood.

\noindent Then for any $y\in B(x_0, t_0, AR(x_0, t_0)^{-1/2})$, we
have
\begin{equation*}
\frac{R(y, t_0)}{R(x_0, t_0)}\leq \Lambda.
\end{equation*}

{\bf Sketch of Proof} One can easily check that the proof of
[BBB$^+$, Theorem 6.1.1] and [BBM, Theorem 6.4] can be adapted to
our situation. For some of the details one can also consult Step 2
of proof of [CZ2, Theorem 4.1] (for the smooth (without surgery)
case) and Step 3 of proof of [CZ2, Proposition 5.4] (for the
surgical case).

\hspace *{0.4cm}

{\bf Remark} For the  estimate above, under a parabolic rescaling of
the metrics,  $c, \varrho, P, R$, etc. will change in general, and
$Q$ will change with the same scaling factor as $R$ does, but
$\Lambda$ is scaling invariant.

\hspace *{0.4cm}

{\bf Acknowledgements } I would like to thank Prof. Sylvain Maillot
for sending me his preprint [BBM] with Bessi$\grave{e}$res and
Besson before it was posted on arXiv. I'm  also very grateful to
Prof. Xi-Ping Zhu for helpful communications, in particular for his
answering my questions on his paper [CaZ] with H.-D. Cao, which is
relevant to Claim 2 in the proof of Proposition 3.6 here.


\hspace *{0.4cm}

{\bf Reference}

\bibliography{1}[BBB$^+$] L. Bessi$\grave{e}$res, G. Besson, M. Boileau, S.
Maillot and J.Porti, Geometrisation of 3-manifolds, Europ. Math.
Soc. 2010.

\bibliography{2}[BBM] L. Bessi$\grave{e}$res, G. Besson and S. Maillot, Ricci flow on
 open 3-manifolds and positive scalar curvature,  arxiv:1001.1458; appeared in Geometry and Topology 15 (2011), 927-975.

\bibliography{4}[CaZ] H.-D. Cao, X.-P. Zhu, A complete proof of the
Poincar$\acute{e}$ and geometrization conjectures- application of
the Hamilton-Perelman theory of the Ricci flow, Asian J. Math. 10
(2006), 165-492.

\bibliography{5}[CTZ] B.-L. Chen, S.-H. Tang and X.-P. Zhu,  Complete classification of
compact four-manifolds with positive isotropic curvature, arXiv:0810.1999.

\bibliography{6}[CZ1] B.-L. Chen, X.-P. Zhu, Uniqueness of the Ricci flow on complete noncompact manifolds, J. Diff. Geom.
74 (2006), 119-154.

\bibliography{7}[CZ2] B.-L. Chen, X.-P. Zhu, Ricci flow with surgery
on  four-manifolds with positive isotropic curvature, J. Diff. Geom.
74 (2006), 177-264.

\bibliography{8}[CCG$^+$08] B. Chow, S.-C. Chu, D. Glickenstein, C. Guenther, J. Isenberg, T. Ivey,
D. Knopf, P. Lu, F. Luo, and  L. Ni, The Ricci flow: techniques and
applications. Part II. Analytic aspects. Mathematical Surveys and
Monographs, 144. American Mathematical Society, Providence, RI,
2008.

\bibliography{9}[CLN] B. Chow, P. Lu and L. Ni, Hamilton's Ricci
flow, Science Press and  AMS 2006.

\bibliography{10}[H1] R. Hamilton, Three-manifolds with positive
Ricci curvature, J. Diff. Geom. 17(1982), 255-306.

\bibliography{11}[H2] R. Hamilton, Four-manifolds with positive
curvature operator, J. Diff. Geom. 24(1986), 153-179.

\bibliography{12}[H3] R. Hamilton, The Harnack estimate for the Ricci flow, J. Diff. Geom. 37 (1993), 225-243.

\bibliography{13}[H4] R. Hamilton, A compactness property for solutions of the Ricci flow, Amer. J.
Math. 117 (1995), no.3, 545-572.

\bibliography{14}[H5] R. Hamilton, Four manifolds with positive isotropic curvature, Comm. Anal. Geom. 5 (1997), 1-92.

\bibliography{15}[Hu1] H. Huang, Complete 4-manifolds with uniformly positive isotropic curvature, arXiv:0912.5405v6.

\bibliography{16}[Hu2] H. Huang, Four orbifolds with positive isotropic curvature, arXiv:1107.1469v2.

\bibliography{17}[LT] P. Lu, G. Tian, Uniqueness of standard
solutions in the work of Perelman,
http://math.berkeley.edu/~lott/ricciflow/StanUniqWork2.pdf

\bibliography{18}[KL] B. Kleiner, J. Lott, Notes on Perelman's
papers, Geom. Topol. 12 (2008), 2587-2855.

\bibliography{19}[MM] M. Micallef, J. D. Moore, Minimal two-spheres
and the topology of manifolds with positive curvature on totally
isotropic two-planes, Ann. Math. (2) 127 (1988), 199-227.

\bibliography{20}[MW] M. Micallef, M. Wang, Metrics with nonnegative isotropic curvature, Duke. Math. J. 72 (1993), no. 3.  649-672.

\bibliography{22}[MT] J. Morgan, G. Tian, Ricci flow and the
Poincar$\acute{e}$ conjecture, Clay Mathematics Monographs 3, Amer.
Math. Soc., 2007.

\bibliography{23}[P1] G. Perelman, The entropy formula for the Ricci flow and its geometric applications,
arXiv:math.DG/0211159.

\bibliography{24}[P2] G. Perelman, Ricci flow with surgery on three-manifolds,

arXiv:math.DG/0303109.

\bibliography{25}[S] W.-X. Shi, Deforming the metric on complete
Riemannian manifolds, J. Diff. Geom. 30 (1989), 223-301.

\bibliography{26} [Z]  Q. S. Zhang, Sobolev inequalities, heat
kernels under Ricci flow, and the Poincar$\acute{e}$ conjecture, CRC
Press 2011.

\vspace *{0.4cm}


\end{document}